\documentclass[12pt]{amsart}

\usepackage{tikz}
\pgfmathsetmacro{\xcoord}{cos(60)}
\pgfmathsetmacro{\ycoord}{sin(60)}
\newcommand{\coords}{\pgftransformcm{1}{0}{\xcoord}{\ycoord}{\pgfpointorigin}}
\newcommand{\hextiling}{
\coords
\draw[fill=\zerocolour] (0,0) -- (20,0) -- (0,20) -- cycle;
\path (0,20) +(0,-0) -- node[left] {$d$} +(0,-1);
\path (0,19) +(0,-1) -- node[left] {$c$} +(0,-5);
\path (0,18) +(0,-5) -- node[left] {$b$} +(0,-8);
\path (0,17) +(0,-8) -- node[left] {$a$} +(0,-11);
\path (20,0) ++(-0,0) +(-0,0) -- node[right] {$h$} +(-1,1);
\path (20,0) ++(-1,1) +(-1,1) -- node[right] {$g$} +(-7,7);
\path (20,0) ++(-2,2) +(-7,7) -- node[right] {$f$} +(-8,8);
\path (20,0) ++(-3,3) +(-8,8) -- node[right] {$e$} +(-11,11);
\path (20,0) +(-0,0) -- node[below] {$\ell$} +(-10,0);
\path (19,0) +(-10,0) -- node[below] {$k$} +(-11,0);
\path (18,0) +(-11,0) -- node[below] {$j$} +(-13,0);
\path (17,0) +(-13,0) -- node[below] {$i$} +(-15,0);
\draw[fill=\rhombuscolour] (9,0) rectangle node {$h$} +(1,1);
\draw[fill=\rhombuscolour] (9,1) ++(0,0) ++(0,1) -- ++(1,-1) -- (19,1) -- ++(-1,1) -- cycle;
\draw[fill=\onecolour] (9,1) -- +(1,0) -- +(0,1) -- cycle;
\draw[fill=\rhombuscolour] (9,1) -- ++(0,1) -- ++(-1,1) -- ++(0,-1) -- cycle;
\draw[fill=\onecolour] (8,2) -- ++(0,1) -- ++(-1,0) -- cycle;
\draw[fill=\rhombuscolour] (7,0) rectangle node {$p$} +(1,2);
\draw[fill=\rhombuscolour] (7,2) ++(0,0) ++(0,1) -- ++(1,-1) -- (8,2) -- ++(-1,1) -- cycle;
\draw[fill=\onecolour] (7,2) -- +(1,0) -- +(0,1) -- cycle;
\draw[fill=\rhombuscolour] (7,2) -- ++(0,1) -- ++(-2,2) -- ++(0,-1) -- cycle;
\draw[fill=\onecolour] (5,4) -- ++(0,1) -- ++(-1,0) -- cycle;
\draw[fill=\rhombuscolour] (4,0) rectangle node {$s$} +(1,3);
\draw[fill=\rhombuscolour] (4,4) ++(0,0) ++(0,1) -- ++(1,-1) -- (5,4) -- ++(-1,1) -- cycle;
\draw[fill=green] (4,3) -- +(1,0) -- +(1,1) -- +(0,2) -- +(-1,2) -- +(-1,1) -- cycle;
\draw[fill=\rhombuscolour] (3,4) -- ++(0,1) -- ++(-0,0) -- ++(0,-1) -- cycle;
\draw[fill=\onecolour] (3,4) -- ++(0,1) -- ++(-1,0) -- cycle;
\draw[fill=\rhombuscolour] (1,0) rectangle node {$u$} +(1,4);
\draw[fill=\rhombuscolour] (1,4) ++(0,0) ++(0,1) -- ++(1,-1) -- (3,4) -- ++(-1,1) -- cycle;
\draw[fill=\onecolour] (1,4) -- +(1,0) -- +(0,1) -- cycle;
\draw[fill=\rhombuscolour] (1,4) -- ++(0,1) -- ++(-1,1) -- ++(0,-1) -- cycle;
\draw[fill=\rhombuscolour] (7,3) rectangle node {$m$} +(1,5);
\draw[fill=\rhombuscolour] (7,8) ++(0,0) ++(0,1) -- ++(1,-1) -- (12,8) -- ++(-1,1) -- cycle;
\draw[fill=\onecolour] (7,8) -- +(1,0) -- +(0,1) -- cycle;
\draw[fill=\rhombuscolour] (7,8) -- ++(0,1) -- ++(-0,0) -- ++(0,-1) -- cycle;
\draw[fill=\onecolour] (7,8) -- ++(0,1) -- ++(-1,0) -- cycle;
\draw[fill=\rhombuscolour] (4,5) rectangle node {$q$} +(1,2);
\draw[fill=\rhombuscolour] (4,8) ++(0,0) ++(0,1) -- ++(1,-1) -- (7,8) -- ++(-1,1) -- cycle;
\draw[fill=green] (4,7) -- +(1,0) -- +(1,1) -- +(0,2) -- +(-1,2) -- +(-1,1) -- cycle;
\draw[fill=\rhombuscolour] (3,8) -- ++(0,1) -- ++(-0,0) -- ++(0,-1) -- cycle;
\draw[fill=\onecolour] (3,8) -- ++(0,1) -- ++(-1,0) -- cycle;
\draw[fill=\rhombuscolour] (2,5) rectangle node {$t$} +(1,2);
\draw[fill=\rhombuscolour] (2,8) ++(0,0) ++(0,1) -- ++(1,-1) -- (3,8) -- ++(-1,1) -- cycle;
\draw[fill=green] (2,7) -- +(1,0) -- +(1,1) -- +(0,2) -- +(-1,2) -- +(-1,1) -- cycle;
\draw[fill=\rhombuscolour] (1,8) -- ++(0,1) -- ++(-1,1) -- ++(0,-1) -- cycle;
\draw[fill=\rhombuscolour] (6,9) rectangle node {$n$} +(1,1);
\draw[fill=\rhombuscolour] (6,10) ++(0,0) ++(0,1) -- ++(1,-1) -- (10,10) -- ++(-1,1) -- cycle;
\draw[fill=\onecolour] (6,10) -- +(1,0) -- +(0,1) -- cycle;
\draw[fill=\rhombuscolour] (6,10) -- ++(0,1) -- ++(-1,1) -- ++(0,-1) -- cycle;
\draw[fill=\onecolour] (5,11) -- ++(0,1) -- ++(-1,0) -- cycle;
\draw[fill=\rhombuscolour] (2,9) rectangle node {$r$} +(1,2);
\draw[fill=\rhombuscolour] (2,11) ++(0,0) ++(0,1) -- ++(1,-1) -- (5,11) -- ++(-1,1) -- cycle;
\draw[fill=\onecolour] (2,11) -- +(1,0) -- +(0,1) -- cycle;
\draw[fill=\rhombuscolour] (2,11) -- ++(0,1) -- ++(-2,2) -- ++(0,-1) -- cycle;
\draw[fill=\rhombuscolour] (4,12) rectangle node {$o$} +(1,2);
\draw[fill=\rhombuscolour] (4,14) ++(0,0) ++(0,1) -- ++(1,-1) -- (6,14) -- ++(-1,1) -- cycle;
\draw[fill=\onecolour] (4,14) -- +(1,0) -- +(0,1) -- cycle;
\draw[fill=\rhombuscolour] (4,14) -- ++(0,1) -- ++(-4,4) -- ++(0,-1) -- cycle;
}

\usepackage{amssymb}
\usepackage{ytableau}
\usepackage[normalem]{ulem} 

\usepackage{comment}
\usepackage{hyperref}

\usepackage{amsthm,amssymb}
\usepackage{epsfig}
\usepackage{accents}
\usepackage[arrow]{xy}

\usepackage{amssymb,latexsym,cite,epsf,graphics}

\usepackage[margin=1.0in]{geometry}
\usepackage{graphicx,color}

\usepackage{fancyhdr}
\fancypagestyle{title}{
\fancyhf{}
\fancyhead[L]{KTiles draft}
\fancyhead[R]{\today}
}

\definecolor{darkred}{rgb}{1,0,0}

\newtheorem{theorem}{Theorem}[section]
\newtheorem{thm}[theorem]{Theorem}

\theoremstyle{remark}

\theoremstyle{definition}

\newtheorem{example}[theorem]{Example}

\newtheorem*{theorem*}{Theorem}

\def\G{{{\mathcal G}}}

\def\cc{\mathbf{c}}

\title{
Puzzles in $K$-homology of Grassmannians
}

\numberwithin{equation}{section}

\renewcommand{\b}{\vec{b}}
\renewcommand{\c}{\mathbf{c}}

\newcommand{\g}{\gamma}
\renewcommand{\G}{\widetilde G}
\renewcommand{\d}{\searrow}
\newcommand{\D}{\Delta}

\renewcommand{\r}{\rho}

\newcommand{\abs}[1]{\left\lvert#1\right\rvert}
\newcommand{\set}[1]{{\left\{#1\right\}}}

\DeclareMathOperator{\shape}{shape}
\DeclareMathOperator{\content}{content}
\DeclareMathOperator{\row}{row}

\def\r{{\mathfrak {row}}}

\newcommand{\figureprefix}{a}
\newcommand{\fig}[1]{Figure~\figref{#1}}
\newcommand{\figref}[1]{\ref{fig:\figureprefix:#1}}
\newcommand{\lbl}[1]{\label{fig:\figureprefix:#1}}
\newcommand{\latin}[1]{\textsl{#1}}

\newcommand{\defn}[1]{\textbf{#1}}

\newcommand{\lmu}[1]{#1_{\la\mu}^\nu}
\renewcommand{\c}{\lmu{c}}
\renewcommand{\d}{\lmu{d}}
\renewcommand{\cc}{\lmu{\tilde c}}
\newcommand{\dd}{\lmu{\tilde d}}
\newcommand{\DD}{\lmu{\Delta}}
\newcommand{\sign}{(-1)^{\abs\nu-\abs\la-\abs\mu}}

\newcommand{\la}{\lambda}

\newcommand*\circled[1]{%
  \begin{tikzpicture}
      \node[draw,circle,inner sep=1pt] {$#1$};
   \end{tikzpicture}}

\newcommand{\tri}{\blacktriangle}
\newcommand{\hex}{\smash{\raisebox{-1pt}{\begin{tikzpicture}[scale=0.2]
\draw[fill=green] (0,0) -- ++(60:1) -- ++(120:1) -- ++(180:1) -- ++(240:1) -- ++(300:1) -- cycle;
\end{tikzpicture}}}}
\renewcommand{\tri}{\smash{\raisebox{-1pt}{\begin{tikzpicture}[scale=0.2]
\draw[fill=orange] (0,0) -- ++(120:2) -- ++(240:2) -- cycle;
\end{tikzpicture}}}}
\newcommand{\K}{\smash{\raisebox{-1pt}{\begin{tikzpicture}[scale=0.2]
\draw[fill=cyan] (0,0) -- ++(60:2) -- ++(180:2) -- cycle;
\end{tikzpicture}}}}
\newcommand{\hexR}{\smash{\raisebox{-1pt}{\begin{tikzpicture}[scale=0.2]
\draw[fill=magenta] (0,0) -- ++(60:1) -- ++(120:1) -- ++(180:1) -- ++(240:1) -- ++(300:1) -- cycle;
\draw[fill=white] (120:1) circle (0.5);
\end{tikzpicture}}}}

\renewcommand{\star}{\smash{\raisebox{0pt}{\begin{tikzpicture}[scale=0.2]
\draw[fill=\onecolour] (0,0) -- ++(60:1) -- ++(180:1) -- cycle;
\end{tikzpicture}}}}

\renewcommand*\circled[1]{\tikz[baseline=(char.base)]{
            \node[shape=circle,draw,inner sep=2pt] (char) {$#1$};}}

\renewcommand{\r}[1]{to node {0} ++(60*#1:1)}

\renewcommand{\b}[1]{to node {1} ++(60*#1:1)}

\newcommand{\bb}[1]{to node {2} ++(60*#1:1)}

\newcommand{\zerocolour}{red!40}
\newcommand{\onecolour}{blue!70}
\newcommand{\rhombuscolour}{yellow!30}
\begin{document}

\author{Pavlo Pylyavskyy}
\address{\hspace{-.3in} Department of Mathematics, University of Minnesota,
Minneapolis, MN 55414, USA}
\email{ppylyavs@umn.edu}

\author{Jed Yang}
\address{\hspace{-.3in} Department of Computer Science, Carleton College,
Northfield, MN 55057, USA}
\email{jyang@carleton.edu}

\date{\today
}




\begin{abstract}
Knutson, Tao, and Woodward~\cite{KTW} formulated a Littlewood--Richardson rule for the cohomology ring of Grassmannians in terms of puzzles.
Vakil~\cite{Vakil} and Wheeler--Zinn-Justin~\cite{WZ} have found additional triangular puzzle pieces that allow one to express structure constants for $K$-theory of Grassmannians.
Here we introduce two other puzzle pieces of hexagonal shape,
each of which gives a Littlewood--Richardson rule for $K$-homology of Grassmannians.
We also explore the corresponding eight versions of $K$-theoretic Littlewood--Richardson tableaux.
\end{abstract}

\ \vspace{-.1in}

\maketitle
\thispagestyle{empty}


\section{Introduction}

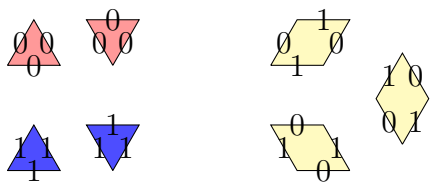
\begin{figure}[hbtp]
\begin{tikzpicture}[scale=0.7]
\draw[fill=\zerocolour] (0,2)
   \r0\r2\r4;
\draw[fill=\zerocolour] (2,2)
   \r1\r3\r5;

\draw[fill=\onecolour] (0,0)
   \b0\b2\b4;
\draw[fill=\onecolour] (2,0)
   \b1\b3\b5;

\begin{scope}[shift={(5,0)}]
\draw[fill=\rhombuscolour] (0,2)
   \b0\r1\b3\r4;
\draw[fill=\rhombuscolour] (1.5,0)
   \b2\r3\b5\r0;
\draw[fill=\rhombuscolour] (2.5,0.5)
   \b1\r2\b4\r5;
\end{scope}
\end{tikzpicture}
\caption{The Knutson--Tao--Woodward tiles.}
\label{LR-tilesintro}
\end{figure}

Cohomology rings of flag varieties are a major object of interest in algebraic geometry, see \cite{Ful, Man} for an exposition. Perhaps the most well-studied and well-understood examples are the cohomology rings of Grassmannians, with a distinguished basis of Schubert classes. A {\it {Littlewood--Richardson rule}} is a combinatorial way to compute the structure constants for this basis. Equivalently, those are the same structure constants $\c$ with which certain symmetric functions -- {\it {Schur functions}} $s_{\lambda}$ -- multiply:
$s_{\lambda} s_{\mu} = \sum_{\nu} \c s_{\nu}.$
In their groundbreaking work  Knutson, Tao, and Woodward \cite{KT, KTW} introduced {\it {puzzles}}, which allow for a powerful formulation of the Littlewood--Richardson rule. Puzzles are tilings of triangular boards with specified boundary conditions by a set of tiles shown in Figure~\ref{LR-tilesintro}. Using puzzles Knutson, Tao, and Woodward studied the faces of the Klyachko cone. 

There is a cohomology theory for each one-dimensional group law, see \cite{Haz, LZ}. For the additive group law $x \oplus y = x+y$ one has the usual cohomology, while the multiplicative group law $x \oplus y = x + y + xy$ gives the {\it {$K$-theory}}.
$K$-theory of Grassmannians was extensively studied, starting with the works of Lascoux and Sch¨utzenberger. In \cite{LS} they introduced the {\it {Grothendieck polynomials}} as representatives of $K$-theory classes of structure sheaves of Schubert varieties.
Fomin and Kirillov \cite{FK} studied those from combinatorial point of view, introducing the {\it {stable Grothendieck polynomials}} $G_{\lambda}$.  Stable Grothendieck polynomials  are symmetric power series that form a rather precise $K$-theoretic analog of Schur functions: their multiplicative structure constants are the same as those for classes of the structure sheaves of Schubert varieties in the corresponding $K$-theory ring.

\begin{figure}[hbtp]
\begin{tikzpicture}[scale=0.7]

\begin{scope}[shift={(-0.5,0)}]
   \node at (0,0.9) {$\K$};
   \draw (0,0) \b1\r1\b3\r3\b5\r5;
\end{scope}

\begin{scope}[shift={(2,0)}]
   \node at (1,0.5) {$\tri$};
   \draw (0,0) \b0\r0\b2\r2\b4\r4;
\end{scope}

\begin{scope}[shift={(7,0)}]
   \node[below] at (120:1) {$\hex$};
   \draw (0,0) \r1\b2\r3\b4\r5\b0;
\end{scope}

\begin{scope}[shift={(10,0)}]
   \node[below] at (60:1) {$\hexR$};
   \draw (0,0) \r0\b1\r2\b3\r4\b5; 
\end{scope}
\end{tikzpicture}
\caption{The four $K$-theoretic tiles.}
\label{K-tilesintro}
\end{figure}
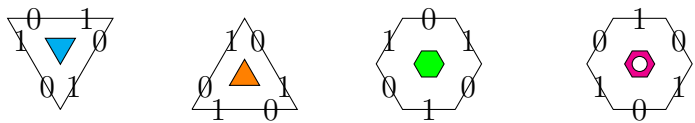

The first $K$-theoretic Littlewood--Richardson rule was obtained by Buch in \cite{Buch}. Vakil \cite{Vakil} has extended puzzles to $K$-theory, giving a puzzle version of the rule. His extension works by adding a single additional tile to the set of tiles from the work of Knutson, Tao and Woodward \cite{KTW}. Later, Wheeler and Zinn-Justin
found an alternative $K$-theoretic tile, that gives the structure constants of dual $K$-theory in an appropriate sense, see \cite{WZ}. Both Vakil and Wheeler-Zinn-Justin tiles have triangular shape and can be seen in Figure~\ref{K-tilesintro}.

In this work we present {\underline {two new tiles}}, adding either one of which  to the standard collection allows to recover structure constants of the Schubert basis in the {\it {$K$-homology}} ring of the Grassmannians, as studied by Lam and Pylyavskyy \cite{LP}.  Equivalently, the corresponding puzzles produce
a combinatorial rule for the {\it {coproduct}} structure constants of the stable Grothendieck polynomials. The first such rule was obtained by Buch in \cite{Buch}. The tiles have {\underline {hexagonal shape}} and can be seen in Figure~\ref{K-tilesintro}.

The paper proceeds as follows. In Section~\ref{sec:1} we recall the known results on the cohomology ring of Grassmannians, including tableaux and puzzles formulations of the Littlewood--Richardson rule. In Section~\ref{sec:2} we recall the $K$-theoretic version of the story, and state our main results regarding the two new hexagonal tiles. We also systematize the eight different tableaux formulations of the $K$-theoretic Littlewood--Richardson rule, some of which are new. The proofs are postponed to Section~\ref{sec:3}. In Section~\ref{sec:4} we conclude with remarks, including the relation of our work to that of Pechenik and Yong \cite{genomic} on genomic tableaux.

\section{Puzzles and tableaux} \label{sec:1}

\subsection{Cohomology of Grassmannians}

A \defn{partition} $\la=(\la_1,\la_2,\dotsc,\la_k)$ is a weakly decreasing sequence $\la_1\geq\la_2\geq\dotsb\geq\la_k\geq 0$
of finitely many nonnegative integers.
The Young diagram, or simply, \defn{diagram}, of $\la$ is
a collection of boxes, top and left justified, with $\la_i$ boxes in Row~$i$.
For example, is the diagram of the partition $\la=(4,3,1)$.
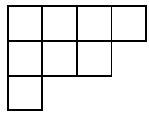
\[\tiny\ydiagram{4,3,1}\]
If $\la$ is a partition whose diagram fits inside that of partition $\nu$,
the \defn{skew} diagram of \defn{shape} $\nu/\la$
is the diagram consisting of the boxes of the diagram of $\nu$ outside that of $\la$.
For example, the following is the diagram of $(4,3,1)/(2,1)$.
\[\tiny\ydiagram{2+2,1+2,1}\]
Given a (possibly skew) diagram and a set $V$,
a \defn{$V$-tableau} $T$ is a filling of the boxes with values in~$V$.
If $V$ is omitted, it is understood that $V$ is the positive integers.
The \defn{shape} of $T$, denoted $\shape(T)$, is the shape of the diagram.
We say that $T$ is \defn{semistandard} if the values are weakly increasing from left to right in rows and strictly increasing from top to bottom in columns.
The \defn{reverse row word} of $T$, denoted $\row(T)$, is the sequence of values of $T$,
read row by row, top to bottom, right to left.
For example,
\[T=\tiny\begin{ytableau}
\none & \none & 1 & 1 \\
\none & 2 & 2 \\
1\\
\end{ytableau}\]
is a semistandard tableau with $\shape(T)=(4,3,1)/(2,1)$ and $\row(T)=11221$.

Let $x_1,x_2,\dotsc$ be commutative variables,
and let $x^T$ denote the monomial $x_{w_1}x_{w_2}\dotsb x_{w_r}$ where $\row(T)=w_1w_2\dotsb w_r$.
The \defn{Schur polynomial} $s_\la$ is given by
\[ s_\la(x)=\sum_T x^T,\]
where the sum runs over all semistandard tableaux $T$ of shape~$\la$.
It is well known that $s_\la$ is symmetric
and $\set{s_\la}_\la$ is a linear basis for the space of all symmetric polynomials
(see \latin{e.g.} \cite{EC2}).
We may therefore expand the product $s_\la s_\mu$ uniquely as a sum of Schur polynomials $s_\nu$ as
\[ s_\la s_\mu = \sum_\nu \c s_\nu. \]
It turns out that $\c$, called the \defn{Littlewood--Richardson coefficient},
is a nonnegative integer,
and is zero whenever $\abs\nu\neq\abs\la+\abs\mu$,
where $\abs\la$ is the number of boxes of $\la$.
In other words, we can let the sum above run over only $\nu$ such that $\abs\nu=\abs\la+\abs\mu$.
This implies that the sum has finitely many terms.

The Littlewood--Richardson coefficients are ubiquitous, appearing naturally in a variety of contexts,
including the study of Schubert calculus
and representation theory of symmetric groups and of general linear groups.
There are also many combinatorial rules for computing~$\c$.
In what follows we recall three rules,
two involving counting tableaux and one involving counting puzzles.

\subsection{Tableau versions of the Littlewood--Richardson rule}

Let $w=w_1w_2\dotsb w_r$ be a sequence of positive integers.
The \defn{content} of $w$, denoted $\content(w)$,
is $(m_1,m_2,\dotsc,m_k)$ such that $m_i$ is the number of occurrences of $i$ in the sequence~$w$.%
\footnote{For example, if $w=\row(T)$, then in the monomial $x^T$, the exponent of $x_i$ is~$m_i$.}
We say $w$ is \defn{ballot} if $\content(w_1\dotsb w_i)$ is a partition for every $i$.
In other words, in every initial segment of $w$, the number $j$ occurs at least as many times as the number $j+1$.
The \defn{content} of $T$, denoted $\content(T)$, is simply $\content(\row(T))$.
We say that $T$ is \defn{ballot} if $\row(T)$ is.

\begin{thm}[Littlewood--Richardson rule, skew version]
For partitions $\la,\mu,\nu$ such that $\abs\nu=\abs\mu+\abs\nu$,
the coefficient $\c$ is the number of semistandard ballot tableaux of shape $\nu/\lambda$ and content $\mu$.
\end{thm}

\begin{example}
Let $\la=(2,1)$, $\mu=(3,2)$, and $\nu=(4,3,1)$ in the following examples.
The following are the (only) two ways to fill according to the Littlewood--Richardson rule.
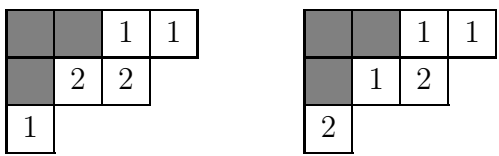
\begin{center}
\begin{ytableau}
*(gray) & *(gray) & 1 & 1 \\
*(gray) & 2 & 2 \\
1\\
\end{ytableau}\hspace{.5in}
\begin{ytableau}
*(gray) & *(gray) & 1 & 1 \\
*(gray) & 1 & 2 \\
2\\
\end{ytableau}
\end{center}
This shows that $\c=2$.
For visual purposes,
we gray out the boxes corresponding to $\la$ instead of removing them.
(This will be useful later when we temporarily write numbers in removed boxes.)
\end{example}

Given two partitions $\la$ and $\mu$,
let the $\oplus$ diagram of shape $\mu\oplus\la$ be obtained by putting the diagrams of $\mu$ and $\la$ corner to corner,
with $\mu$ to the lower left and $\la$ to the upper right.
For example,
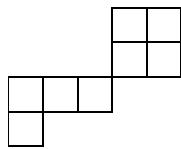
\[ \tiny\ydiagram{3+2,3+2,3,1} \]
is a diagram of shape $(3,1)\oplus (2,2)$.

\begin{thm}[Littlewood--Richardson rule, $\oplus$ version]
For partitions $\la,\mu,\nu$ such that $\abs\nu=\abs\mu+\abs\nu$,
the coefficient $\c$ is the number of semistandard ballot tableaux of shape $\mu\oplus\la$ and content $\nu$.
\end{thm}

\begin{example}
We continue with $\la,\mu,\nu$ from the example above.
The following are the two corresponding fillings using the $\oplus$ version of the Littlewood--Richardson rule.
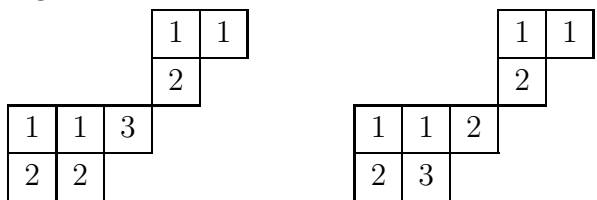
\begin{center}
\begin{ytableau}
\none & \none & \none & 1 & 1 \\
\none & \none & \none & 2 \\
1 & 1 & 3\\
2 & 2
\end{ytableau}\hspace{.5in}
\begin{ytableau}
\none & \none & \none & 1 & 1 \\
\none & \none & \none & 2 \\
1 & 1 & 2\\
2 & 3
\end{ytableau}
\end{center}
These are displayed in the same order under the bijection that is described in later sections.
\end{example}

Of course, any $\oplus$ diagram $\mu\oplus\la$ is also a skew diagram of shape
\[(\la_1+\mu_1,\dotsc,\la_k+\mu_1,\mu_1,\dotsc,\mu_k).\]
Nevertheless, we think of these classes of shapes separately,
since we will have pairs of tableaux rules, one involving shape $\nu/\la$ and one involving shape $\mu\oplus\la$.
We refer to $\nu/\la$ as skew shape (and use grayed out boxes)
and refer to $\mu\oplus\la$ as $\oplus$ shape (without using grayed out boxes).

\subsection{Puzzle version of the Littlewood--Richardson rule}
Let $n\geq k$ be positive integers.
Refer to the partition of $k$ rows of length $n-k$ as the \defn{ambient rectangle}.
From now on, we consider only partitions whose diagrams fit inside this ambient rectangle.
(To consider bigger partitions, simply specify a larger ambient rectangle.)
On the lower right boundary of a partition inside the ambient rectangle,
write
a $0$ on each horizontal edge and a $1$ on each vertical edge (see \fig{binary string}). 
A binary string of length $n$ with $k$ ones and $n-k$ zeros
is obtained by reading these numbers from top right to bottom left.

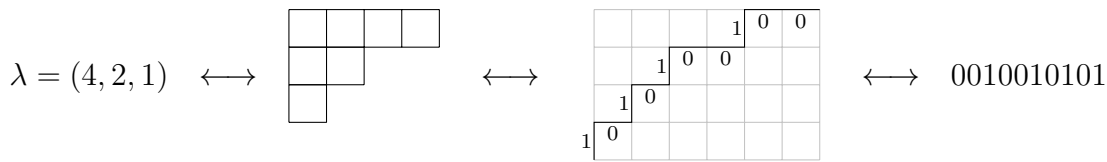
\begin{figure}[hbtp]
\newcommand{\bij}{\quad\ensuremath{\longleftrightarrow}\quad}
$\la=(4,2,1)$%
\bij
\begin{tikzpicture}[scale=0.5,baseline={(0,-1)}]
\draw (0,0) grid ++(4,-1);
\draw (0,-1) grid ++(2,-1);
\draw (0,-2) grid ++(1,-1);
\end{tikzpicture}
%
\bij
\begin{tikzpicture}[scale=0.5,baseline={(0,-1)}]
\renewcommand{\r}{to node[shift={(0,-0.15)}] {\tiny 0} ++(180:1)}
\renewcommand{\b}{to node[shift={(-0.1,0)}] {\tiny 1} ++(270:1)}
\draw[color=black!23] (0,0) grid (-6,-4);
\draw (0,0) \r\r\b\r\r\b\r\b\r\b;
\end{tikzpicture}
\bij%
{%
   \renewcommand{\r}{0}%
   \renewcommand{\b}{1}%
   \r\r\b\r\r\b\r\b\r\b
}
\caption{Bijection between partitions, Young diagrams, and binary strings; $n=10$, $k=4$.}
\lbl{binary string}
\end{figure}
Here we consider tilings on the triangular lattice.
Knutson, Tao, and Woodward \cite{KTW} introduced the following \defn{puzzle pieces} (see \fig{LR-tiles}).

\begin{itemize}
\item $0$-triangle: unit triangle with edges labelled by $0$, two rotations;
\item $1$-triangle: unit triangle with edges labelled by $1$, two rotations; and
\item rhombus: formed by gluing two adjacent unit triangles together,
with edges labelled by $0$ if clockwise of an acute angle and $1$ if clockwise of an obtuse angle, three rotations.
\end{itemize}

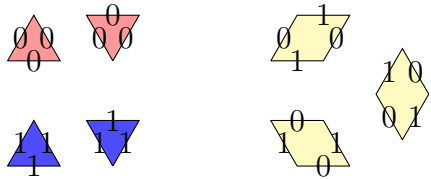
\begin{figure}[hbtp]
\begin{tikzpicture}[scale=0.7]
\draw[fill=\zerocolour] (0,2)
   \r0\r2\r4;
\draw[fill=\zerocolour] (2,2)
   \r1\r3\r5;

\draw[fill=\onecolour] (0,0)
   \b0\b2\b4;
\draw[fill=\onecolour] (2,0)
   \b1\b3\b5;

\begin{scope}[shift={(5,0)}]
\draw[fill=\rhombuscolour] (0,2)
   \b0\r1\b3\r4;
\draw[fill=\rhombuscolour] (1.5,0)
   \b2\r3\b5\r0;
\draw[fill=\rhombuscolour] (2.5,0.5)
   \b1\r2\b4\r5;
\end{scope}
\end{tikzpicture}
\caption{Puzzle pieces.}
\lbl{LR-tiles}
\end{figure}

A \defn{tiling} is an assembly of (lattice) translated copies of tiles, where edge labels of adjacent tiles must match.
We are interested in tiling an upright triangular region $\DD$ whose
boundary labels of the left, right, and bottom sides, read left-to-right,
are the binary strings corresponding to $\la$, $\mu$, and $\nu$ (see \fig{bdry}).

\begin{figure}[hbtp]
\begin{tikzpicture}[scale=0.5]
\draw (0,0)
   \r0\r0\b0\r0\b0\r0\r0\b0\r0
   \b2\r2\r2\b2\r2\b2\r2\r2\r2
   \b4\r4\b4\r4\b4\r4\r4\r4\r4;
\path (0,0) -- node[below=10] {$\nu$} ++(0:9) -- node[right=10] {$\mu$} ++(120:9) -- node[left=10] {$\la$} cycle;
\end{tikzpicture}
\caption{Boundary $\DD$ with $\la=(2,1,0)$, $\mu=(3,2,0)$, and $\nu=(4,3,1)$.}
\lbl{bdry}
\end{figure}

Littlewood--Richardson coefficients can be calculated by counting puzzle tilings:
\begin{thm}[Knutson--Tao--Woodward \cite{KTW}]
Suppose $\la,\mu,\nu$ are partitions fitting inside an $(n-k)\times k$ ambient rectangle,
with $\abs\nu=\abs\la+\abs\mu$.
The number of puzzle tilings with boundary $\DD$ is $\c$.
\end{thm}

\begin{example}
Continuing with the running example from the previous section, since $\c=2$, there are two tilings of $\DD$:
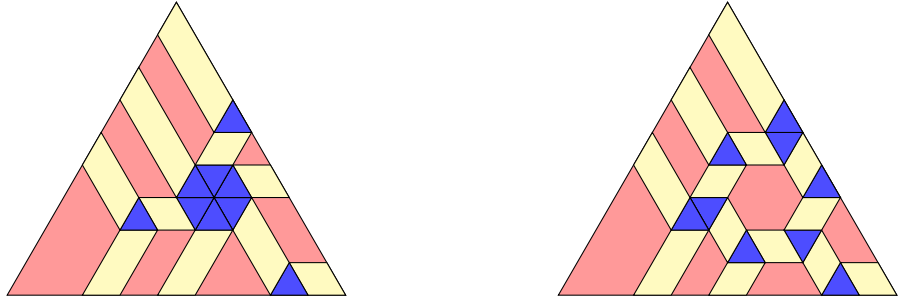
\begin{center}
\begin{tikzpicture}[scale=0.5]
\coords
\draw[fill=\zerocolour] (0,0) -- (9,0) -- (0,9) -- cycle;
\draw[fill=\rhombuscolour] (7,0) rectangle +(1,0);
\draw[fill=\rhombuscolour] (7,0) ++(0,0) ++(0,1) -- ++(1,-1) -- (9,0) -- ++(-1,1) -- cycle;
\draw[fill=\onecolour] (7,0) -- +(1,0) -- +(0,1) -- cycle;
\draw[fill=\rhombuscolour] (7,0) -- ++(0,1) -- ++(-2,2) -- ++(0,-1) -- cycle;
\draw[fill=\onecolour] (5,2) -- ++(0,1) -- ++(-1,0) -- cycle;
\draw[fill=\rhombuscolour] (4,0) rectangle +(1,2);
\draw[fill=\rhombuscolour] (4,2) ++(0,0) ++(0,1) -- ++(1,-1) -- (5,2) -- ++(-1,1) -- cycle;
\draw[fill=\onecolour] (4,2) -- +(1,0) -- +(0,1) -- cycle;
\draw[fill=\rhombuscolour] (4,2) -- ++(0,1) -- ++(-0,0) -- ++(0,-1) -- cycle;
\draw[fill=\onecolour] (4,2) -- ++(0,1) -- ++(-1,0) -- cycle;
\draw[fill=\rhombuscolour] (2,0) rectangle +(1,2);
\draw[fill=\rhombuscolour] (2,2) ++(0,0) ++(0,1) -- ++(1,-1) -- (4,2) -- ++(-1,1) -- cycle;
\draw[fill=\onecolour] (2,2) -- +(1,0) -- +(0,1) -- cycle;
\draw[fill=\rhombuscolour] (2,2) -- ++(0,1) -- ++(-2,2) -- ++(0,-1) -- cycle;
\draw[fill=\rhombuscolour] (4,3) rectangle +(1,0);
\draw[fill=\rhombuscolour] (4,3) ++(0,0) ++(0,1) -- ++(1,-1) -- (6,3) -- ++(-1,1) -- cycle;
\draw[fill=\onecolour] (4,3) -- +(1,0) -- +(0,1) -- cycle;
\draw[fill=\rhombuscolour] (4,3) -- ++(0,1) -- ++(-0,0) -- ++(0,-1) -- cycle;
\draw[fill=\onecolour] (4,3) -- ++(0,1) -- ++(-1,0) -- cycle;
\draw[fill=\rhombuscolour] (3,3) rectangle +(1,0);
\draw[fill=\rhombuscolour] (3,3) ++(0,0) ++(0,1) -- ++(1,-1) -- (4,3) -- ++(-1,1) -- cycle;
\draw[fill=\onecolour] (3,3) -- +(1,0) -- +(0,1) -- cycle;
\draw[fill=\rhombuscolour] (3,3) -- ++(0,1) -- ++(-3,3) -- ++(0,-1) -- cycle;
\draw[fill=\rhombuscolour] (3,4) rectangle +(1,1);
\draw[fill=\rhombuscolour] (3,5) ++(0,0) ++(0,1) -- ++(1,-1) -- (4,5) -- ++(-1,1) -- cycle;
\draw[fill=\onecolour] (3,5) -- +(1,0) -- +(0,1) -- cycle;
\draw[fill=\rhombuscolour] (3,5) -- ++(0,1) -- ++(-3,3) -- ++(0,-1) -- cycle;
\end{tikzpicture}
\hspace{1in}
\begin{tikzpicture}[scale=0.5]
\coords
\draw[fill=\zerocolour] (0,0) -- (9,0) -- (0,9) -- cycle;
\draw[fill=\rhombuscolour] (7,0) rectangle +(1,0);
\draw[fill=\rhombuscolour] (7,0) ++(0,0) ++(0,1) -- ++(1,-1) -- (9,0) -- ++(-1,1) -- cycle;
\draw[fill=\onecolour] (7,0) -- +(1,0) -- +(0,1) -- cycle;
\draw[fill=\rhombuscolour] (7,0) -- ++(0,1) -- ++(-1,1) -- ++(0,-1) -- cycle;
\draw[fill=\onecolour] (6,1) -- ++(0,1) -- ++(-1,0) -- cycle;
\draw[fill=\rhombuscolour] (4,0) rectangle +(1,1);
\draw[fill=\rhombuscolour] (4,1) ++(0,0) ++(0,1) -- ++(1,-1) -- (6,1) -- ++(-1,1) -- cycle;
\draw[fill=\onecolour] (4,1) -- +(1,0) -- +(0,1) -- cycle;
\draw[fill=\rhombuscolour] (4,1) -- ++(0,1) -- ++(-1,1) -- ++(0,-1) -- cycle;
\draw[fill=\onecolour] (3,2) -- ++(0,1) -- ++(-1,0) -- cycle;
\draw[fill=\rhombuscolour] (2,0) rectangle +(1,2);
\draw[fill=\rhombuscolour] (2,2) ++(0,0) ++(0,1) -- ++(1,-1) -- (3,2) -- ++(-1,1) -- cycle;
\draw[fill=\onecolour] (2,2) -- +(1,0) -- +(0,1) -- cycle;
\draw[fill=\rhombuscolour] (2,2) -- ++(0,1) -- ++(-2,2) -- ++(0,-1) -- cycle;
\draw[fill=\rhombuscolour] (5,2) rectangle +(1,1);
\draw[fill=\rhombuscolour] (5,3) ++(0,0) ++(0,1) -- ++(1,-1) -- (6,3) -- ++(-1,1) -- cycle;
\draw[fill=\onecolour] (5,3) -- +(1,0) -- +(0,1) -- cycle;
\draw[fill=\rhombuscolour] (5,3) -- ++(0,1) -- ++(-1,1) -- ++(0,-1) -- cycle;
\draw[fill=\onecolour] (4,4) -- ++(0,1) -- ++(-1,0) -- cycle;
\draw[fill=\rhombuscolour] (2,3) rectangle +(1,1);
\draw[fill=\rhombuscolour] (2,4) ++(0,0) ++(0,1) -- ++(1,-1) -- (4,4) -- ++(-1,1) -- cycle;
\draw[fill=\onecolour] (2,4) -- +(1,0) -- +(0,1) -- cycle;
\draw[fill=\rhombuscolour] (2,4) -- ++(0,1) -- ++(-2,2) -- ++(0,-1) -- cycle;
\draw[fill=\rhombuscolour] (3,5) rectangle +(1,0);
\draw[fill=\rhombuscolour] (3,5) ++(0,0) ++(0,1) -- ++(1,-1) -- (4,5) -- ++(-1,1) -- cycle;
\draw[fill=\onecolour] (3,5) -- +(1,0) -- +(0,1) -- cycle;
\draw[fill=\rhombuscolour] (3,5) -- ++(0,1) -- ++(-3,3) -- ++(0,-1) -- cycle;
\end{tikzpicture}
\end{center}
Here and subsequently, some edges (namely, the edges within a region of $0$-triangles and the $1$-edges of a sequence of rhombi) are omitted to suggest the structure of puzzle tilings.
\end{example}

\section{$K$-theoretic puzzles and tableaux} \label{sec:2}
In this section, we discuss four $K$-theoretic analogues of the Littlewood--Richardson coefficients.
These coefficients can be calculated using four puzzle rules and eight tableaux rules.

\subsection{$K$-theory and $K$-homology of Grassmannians} \label{ss:K-coefficients}

The $K$-theoretic analogue of a Schur function $s_\la$ is the \defn{single stable Grothendieck polynomial} $G_\la$
given by the formula
\[ G_\la = \sum_T(-1)^{\abs T-\abs\la}x^T, \]
where the sum runs over all semistandard set-valued tableaux $T$ of shape $\la$,
and $\abs T$ is the length of $\row(T)$.
The equivalence of this definition to other definitions is established by Buch~\cite{Buch}.

Buch also showed that the linear span of $\set{G_\la}_\la$
is a bialgebra,
with product given by
\[ G_\la G_\mu = \sum_\nu \sign\c G_\nu \]
and coproduct $\D$ given by
\[ \D(G_\nu) = \sum_{\la,\mu}\sign\d G_\la\otimes G_\mu. \]
It turns out that $\c=0$ when $\abs\nu<\abs\la+\abs\mu$
and $\d=0$ when $\abs\nu>\abs\la+\abs\mu$.
So we might as well restrict the first and second sums to the cases where
$\abs\nu\geq\abs\la+\abs\mu$ and
$\abs\nu\leq\abs\la+\abs\mu$, respectively.
Unlike the classical case, this does not immediately show that the sums are finite,
but indeed they are (Corollaries 5.5 and 6.7 of \cite{Buch}).

When $\abs\nu=\abs\la+\abs\mu$,
the number $\c$ is indeed the classical Littlewood--Richardson coefficient described in previous sections.
Since this is the only case where the classical $\c$ is possibly nonzero,
by an abuse of notation, we use the same symbol to denote both.
It is therefore paramount to require $\abs\nu=\abs\la+\abs\mu$ when discussing $\c$ in the classical case.

The following slight variants of $\c$ and $\d$ arise naturally in the study of puzzles.
Let $\G_\la=G_\la\cdot (1-G_1)$.
Define $\cc$ as the unique numbers such that
\[ \G_\la\cdot\G_\mu=\sum_\nu\sign\cc\G_\nu.\]
We again restrict to $\abs\nu\geq\abs\la+\abs\mu$, the only time when $\cc$ is possibly nonzero.

%
%

Finally, let $\dd$ be given by $d_{\la'\mu'}^{\nu'}$, where $\la'$ is the transpose of $\la$,
\latin{i.e.}, mirror the diagram of $\la$ across the line $x+y=0$.
Since the number of boxes is preserved, the only time $\dd$ is possibly nonzero is
when $\abs\nu\leq\abs\la+\abs\mu$.

\subsection{The four $K$-theoretic puzzles} \label{ss:K-puzzles}
Consider the puzzle pieces shown in \fig{K-tiles}.

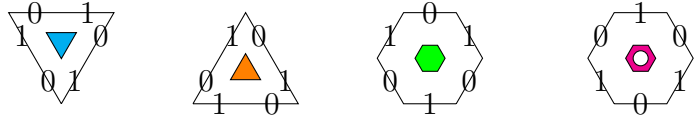
\begin{figure}[hbtp]
\begin{tikzpicture}[scale=0.7]
\begin{scope}[shift={(-0.5,0)}]
   \node at (0,0.9) {$\K$};
   \draw (0,0) \b1\r1\b3\r3\b5\r5;
\end{scope}

\begin{scope}[shift={(2,0)}]
   \node at (1,0.5) {$\tri$};
   \draw (0,0) \b0\r0\b2\r2\b4\r4;
\end{scope}

\begin{scope}[shift={(7,0)}]
   \node[below] at (120:1) {$\hex$};
   \draw (0,0) \r1\b2\r3\b4\r5\b0;
\end{scope}

\begin{scope}[shift={(10,0)}]
   \node[below] at (60:1) {$\hexR$};
   \draw (0,0) \r0\b1\r2\b3\r4\b5; 
\end{scope}
\end{tikzpicture}
\caption{Four additional puzzle pieces.}
\lbl{K-tiles}
\end{figure}

We refer to these puzzle pieces using the corresponding pictograms shown in the figure.
If $X$ is (the pictogram of) an additional puzzle piece, an $X$-puzzle
is a puzzle tiling where, in additional to the usual puzzle pieces, translated copies of $X$ can be used.
There are known interpretations of $\K$-puzzles and $\tri$-puzzles.

\begin{thm}[Vakil \cite{Vakil}] \label{thm:K-tiling}
Suppose $\la,\mu,\nu$ are partitions fitting inside an $(n-k)\times k$ ambient rectangle,
with $\abs\nu\geq\abs\la+\abs\mu$.
The number of $\K$-puzzle tilings with boundary $\DD$ is~$\c$. 
\end{thm}

\begin{thm}[Wheeler--Zinn-Justin \cite{WZ}] \label{thm:tri-tiling}
Suppose $\la,\mu,\nu$ are partitions fitting inside an $(n-k)\times k$ ambient rectangle,
with $\abs\nu\geq\abs\la+\abs\mu$.
The number of $\tri$-puzzle tilings with boundary $\DD$ is~$\cc$. 
\end{thm}

We establish interpretations of $\hex$-puzzles and $\hexR$-puzzles.
\begin{thm} \label{thm:hex-tiling}
Suppose $\la,\mu,\nu$ are partitions fitting inside an $(n-k-1)\times k$ ambient rectangle,%
\footnote{For technical reasons, we require partitions to be slightly smaller.  See Section~\ref{ss:protrude}.}
with $\abs\nu\leq\abs\la+\abs\mu$.
The number of $\hex$-puzzle tilings with boundary $\DD$ is~$\d$. 
\end{thm}

\begin{thm} \label{thm:hexR-tiling}
Suppose $\la,\mu,\nu$ are partitions fitting inside an $(n-k)\times(k-1)$ ambient rectangle,
with $\abs\nu\leq\abs\la+\abs\mu$.
The number of $\hexR$-puzzle tilings with boundary $\DD$ is~$\dd$. 
\end{thm}

\subsection{The eightfold way} \label{ss:K-tableaux}

Like the classical case, where the puzzle rule corresponds to a pair of tableau rules
(involving diagrams of shapes $\nu/\la$ and $\mu\oplus\la$, respectively),
we describe four pairs of $K$-tableau rules corresponding to the four $K$-puzzle rules.

A \defn{set-valued tableau} is a $V$-tableau
where $V$ consists of non-empty subsets of $\set{1,\dotsc,k}$.
To understand the semistandard condition in this context, we agree that
for $A,B\in V$, $A$ is (strictly) less than $B$ if $\max A$ is (strictly) less than $\min B$.
When forming the reverse row word,
a value $A\in V$ is expanded as the numbers in the set $A$, written from largest to smallest.

Buch~\cite{Buch} gives a combinatorial rule for calculating
the $K$-theory Littlewood--Richardson coefficient $\c$
by counting certain set-valued tableaux of $\oplus$ shape.
\begin{thm}[$\K$ rule, $\oplus$ version] \label{thm:K-tab-oplus}
The coefficient $\c$ is the number of semistandard ballot set-valued tableaux of shape $\mu\oplus\la$ and content $\nu$.
\end{thm}

To describe the skew version of the $K$-theory rule, we consider a new kind of tableaux.
A \defn{circle tableau} $T$ is a $V$-tableau
where $V$ consists of $\set{1,\dotsc,k}$ and the \defn{circled} numbers $\{\circled 1,\dotsc,\circled k\}$.

We say $T$ is a \defn{right} (resp., \defn{left}) circle tableau
if each $\circled i$ is the rightmost (resp., leftmost) $i$ or $\circled i$ in its row.
(In other words, for each $i$, only the rightmost (resp., leftmost) $i$ in a row is optionally circled.)
Moreover, circled values may only occur in the bottom $k$ rows (that is, anywhere in shape $\nu/\la$,
bottom half in shape $\mu\oplus\la$).

We say $T$ is \defn{semistandard} if it is semistandard when the circled values are treated as if they are not circled.
Its \defn{content} is $\content(w)$ where $w$ is $\row(T)$ with the circled values omitted.

Let $w$ be an initial segment of $\row(T)$.
If $w$ ends with $\circled i$, replace it with an uncircled $i+1$.
Remove all other circled entries.
Call the result the \defn{incremented erasure} of~$w$.
Analogously, call the result the \defn{unincremented erasure} of~$w$ if the final $\circled i$ is replaced with an uncircled~$i$ instead.
We say that a right (left) circle tableau is \defn{ballot} if all its incremented (unincremented) erasures are ballot.

Pechenik and Yong~\cite{genomic} gives a combinatorial rule for calculating
the $K$-theory Littlewood--Richardson coefficient $\c$
by counting certain \emph{genomic} tableaux of skew shape.
We give an equivalent formulation (see Section~\ref{ss:genomic}) here in terms of circle tableaux.
\begin{thm}[$\K$ rule, skew version] \label{thm:K-tab-skew}
The coefficient $\c$ is the number of semistandard ballot right circle tableaux of shape $\nu/\la$ and content $\mu$.
\end{thm}

An \defn{outer corner} of (the diagram of) a partition $\mu$ is a box whose addition results in a diagram of a partition.

\begin{thm}[$\tri$ rule, $\oplus$ version] \label{thm:tri-tab-oplus}
The coefficient $\cc$ is the number of semistandard ballot set-valued tableaux of shape $\mu^+\oplus\la$ and content $\nu$, where $\mu^+$ is $\mu$ with some number (possibly zero) of its outer corners added.
\end{thm}

\begin{thm}[$\tri$ rule, skew version] \label{thm:tri-tab-skew}
The coefficient $\cc$ is the number of semistandard ballot left circle tableaux of shape $\nu/\la$ and content $\mu$.
\end{thm}

Recall that a circle tableau of shape $\mu\oplus\la$ do not have circles in the rows corresponding to~$\la$.

\begin{thm}[$\hex$ rule, $\oplus$ version] \label{thm:hex-tab-oplus}
The coefficient $\d$ is the number of semistandard ballot right circle tableaux of shape $\mu\oplus\la$ and content $\nu$.
\end{thm}

An \defn{inner corner} of (the diagram of) a partition $\la$ is a box whose removal results in a diagram of a partition.

\begin{thm}[$\hex$ rule, skew version] \label{thm:hex-tab-skew}
The coefficient $\d$ is the number of semistandard ballot set-valued tableaux of shape $\nu/\la^-$ and content $\mu$, where $\la^-$ is $\la$ with some number (possibly zero) of its inner corners removed.
\end{thm}

A circle tableau of shape $\mu\oplus\la$ is \defn{limited}
if it has no $\circled i$ in Row~$i$ of the bottom half for any~$i$.

\begin{thm}[$\hexR$ rule, $\oplus$ version] \label{thm:hexR-tab-oplus}
The coefficient $\dd$ is the number of limited semistandard ballot left circle tableaux of shape $\mu\oplus\la$ and content $\nu$.
\end{thm}

\begin{thm}[$\hexR$ rule, skew version] \label{thm:hexR-tab-skew}
The coefficient $\dd$ is the number of semistandard ballot set-valued tableaux of shape $\nu/\la$ and content $\mu$.
\end{thm}

\section{Proofs} \label{sec:3}

\subsection{Proof of Theorem~\ref{thm:hex-tab-skew}}
Given a sequence $w=(w_1,\dotsc,w_r)$ and an interval $[a,b]$,
let $w|_{[a,b]}$ be the sequence obtained by 
shifting the numbers down to the interval $[1,b-a+1]$ by subtracting $a-1$ from each number $w_i$ in the range $[a,b]$
(and omitting numbers that are out of the range).

\begin{thm}[Buch~\cite{Buch}] \label{thm:buch}
The coefficient $\d$ is the number of semistandard set-valued tableaux $T$ of shape $\nu$
with content $(\la,\mu)=(\la_1,\la_2,\dotsc,\la_k,\mu_1,\dotsc,\mu_k)$,
such that $\row(T)|_{[1,k]}$ 
and $\row(T)|_{[k+1,2k]}$ 
are both ballot.
\end{thm}

For notational convenience, local to this proof only,
a \defn{Buch tableau} is one described in Theorem~\ref{thm:buch}.
and a \defn{PY tableau} is one described in Theorem~\ref{thm:hex-tab-skew}.
There is a simple bijection between Buch tableaux and PY tableaux.

Indeed, let $T$ be a PY tableau.
Increase each number in $T$ by~$k$.
Extend the shape of $T$ to $\nu$ by filling in the first $\la_i$ boxes of $T$ with $i$ in Row~$i$.
The result is clearly a Buch tableau.

Conversely, let $T$ be a Buch tableau.
It is easy to see that,
as $T$ is semistandard and $\row(T)|_{[1,k]}$ is ballot,
the $\la_i$ occurrences of $i$ are exactly in the first $\la_i$ boxes of Row~$i$.
Remove these ``\emph{small}'' numbers.
A remaining ``\emph{big}'' number in Row~$i$ cannot be in the first $\la_i-1$ boxes,
since the $\la_i$-th box contained a small number.
It can be in the $\la_i$-th box only if the $\la_i$-th box in the next row did not contain a small number.
In other words, only if this box is an inner corner of $\la$.
We therefore conclude that the shape of the remaining tableau is $\nu/\la$
with some (possibly zero) inner corners of $\la$ added.
Decrease $k$ from all the remaining numbers to obtain a PY tableau.

This concludes the proof of Theorem~\ref{thm:hex-tab-skew}.

\subsection{Proof of Theorem~\ref{thm:hex-tiling}}

We prove Theorem~\ref{thm:hex-tiling} by establishing a bijection
between $\hex$-puzzles and the tableaux described in Theorem~\ref{thm:hex-tab-skew}.
For notational convenience, we do so by considering an example when $k=4$.
The general case is similar.

\begin{figure}[hbtp]
\begin{tikzpicture}[scale=0.4]
\hextiling
\end{tikzpicture}
\caption{An example tiling.}
\lbl{hex-tiling}
\end{figure}
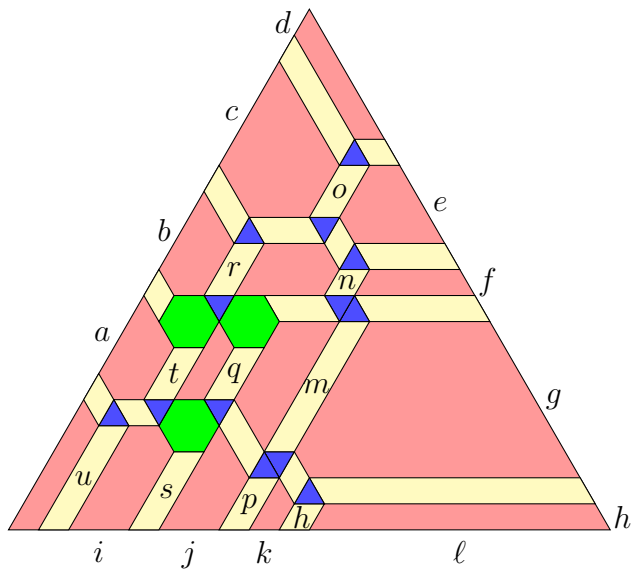

It is easy to see that a $\hex$-puzzle has the structure as in \fig{hex-tiling}.
The sequence of rhombi with label $x$ is called the \defn{$x$-beam}.
By abuse of notation, we also let $x$ denote the length (that is, the number of rhombi) of the $x$-beam.
For each $x$-beam, set $x'$ to $x$.
Increment $x'$ by one if the $x$-beam is capped with a $\hex$ on top (as opposed to a triangle).
In the example above, $t'$, $q'$, and $s'$ are the ones that are incremented.
The boundary also have some length labels.
We use the same labels as those in Tao's ``proof without words'' (see \cite{Vakil}).

We now describe a bijection from the $\hex$-puzzles to skew tableaux.
In the boxes of a diagram of shape $\nu/\la$, fill out according to the following schematic plan
\begin{center}
\begin{ytableau}
*(gray) & *(gray) & *(gray) & *(gray) & 1u\\
*(gray) & *(gray) & *(gray) & 1t & 2s\\
*(gray) & *(gray) & 1r & 2q & 3p\\
*(gray) & 1o & 2n & 3m & 4h\\
\end{ytableau}
\end{center}
where, 
a number $x$ followed by a letter $y$ in the schematic plan means to fill the number $x$ in $y$ consecutive boxes.
If $y'>y$, write an additional $x$ in the \emph{previous} box, without using space.
Circle such a number for easy reference.
The grayed out boxes correspond to $\la$ and may contain circled numbers;
the white boxes correspond to $\nu/\la$ and each has exactly one uncircled number.
Call this tableau~$T$.

\begin{example}
Applying the bijection described to the puzzle results in the following tableau.
\begin{center}
\begin{ytableau}
*(gray) & *(gray) & *(gray) & *(gray) & *(gray) & *(gray) & *(gray) & *(gray) & *(gray) & *(gray) & *(gray) & 1 & 1 & 1 & 1\\
*(gray) & *(gray) & *(gray) & *(gray) & *(gray) & *(gray) & *(gray) & *(gray) \tiny\circled1 & 1 & \text{\tiny$1 \circled2$} & 2 & 2 & 2\\
*(gray) & *(gray) & *(gray) & *(gray) & *(gray) & 1 & \text{\tiny$1 \circled2$} & 2 & 2 & 3 & 3\\
*(gray) & 1 & 1 & 2 & 3 & 3 & 3 & 3 & 3 & 4\\
\end{ytableau}
\end{center}
\end{example}

From the tiling, one could read off certain equalities and inequalities (see \fig{hex-inequalities}).

\textbf{Shape.}
The top left picture shows that $\nu_2+3=1+j+1+k+1+\ell=s+t+1+b+1+c+1+d=s+t+\la_2+3$, or $\nu_2-\la_2=s+t$.
This means that the $s+t$ uncircled numbers we fill in Row 2 of $\nu/\la$
precisely takes up the $\nu_2-\la_2$ boxes.
In other words, the shape is unaffected by the $\hex$ tiles,
except for the possibility of writing $\circled1$ in the shaded boxes, discussed below.

\textbf{Content.}
The top right picture shows that
$s'+1+q'+1+n'=h+1+g+1+h=\mu_2+2$, or $\mu_2=s'+q'+n'$,
leading to $\content(T)=\mu$ where $\circled i$ is treated as~$i$.

\textbf{Ballot.}
The lower left picture shows that
$u'+t'\geq s'+q'\geq p'+m'$.
This directly translates to the ballot condition of~$T$, again by treating $\circled i$ as~$i$.

\textbf{Semistandard.}
The lower right picture shows a final type of inequalities,
which are slightly more complicated.
Let $x\geq_z y$ be a shorthand for $x\geq y+z'-z$.
In other words, $x\geq_z y$ means $x\geq y$ if $z'=z$, and means $x>y$ if $z'=z+1$.
If there are no $\hex$ tiles, the two thick lines in the picture must not cross,
yielding inequalities $b\geq r$ and $b+t\geq r+q$.
Because of the $\hex$ tiles, these inequalities must be strict.
Therefore we get $b\geq_t r$ and $b+t\geq_s r+q$ instead.
These inequalities translate to the semistandard condition of $T$ by considering all pairs of numbers in adjacent boxes.
Also note that if $b=0$, then $b\geq_tr$ says that $t=t'$ (and $r=0$),
so there cannot be a $\circled1$ in Row~$2$ if $\la_2=\la_3$.
Similarly, there cannot be a $\circled1$ in Row~$4$ if $\la_4=0$.
In general, $\circled1$ can only be written in the boxes corresponding to the \emph{inner corners} of~$\la$.

\begin{figure}[hbtp]
\newcommand{\n}{node {$\bullet$}}
\newcommand{\scale}{0.38}
\begin{tikzpicture}[scale=\scale]
\hextiling
\draw[ultra thick] (4,0) \n -- ++(0,3) -- ++(-2,2) -- ++(0,2) -- ++(-2,2) -- ++(0,11) \n;
\draw[ultra thick] (4,0) -- ++(16,0) \n;
\end{tikzpicture}
\begin{tikzpicture}[scale=\scale]
\hextiling
\draw[ultra thick] (5,0) \n -- ++(0,8) -- ++(2,0) -- ++(0,2) -- ++(3,0) \n -- ++(10,-10) \n;
\end{tikzpicture}
\begin{tikzpicture}[scale=\scale]
\hextiling
\draw[ultra thick] (2,4) \n -- ++(3,0) \n -- ++(2,-2) -- ++(1,0) \n -- ++(1,-1) -- ++(1,0) \n;
\draw[ultra thick] (3,8) \n -- ++(2,0) \n -- ++(3,0) \n;
\draw[ultra thick] (3,11) \n -- ++(2,0) -- ++(1,-1) -- ++(1,0) \n;
\end{tikzpicture}
\begin{tikzpicture}[scale=\scale]
\hextiling
\node at (0,13) {$\bullet$};
\draw[ultra thick,->] (0,13) -- ++(0,-3) -- ++(3,-3) -- ++(0,-2) -- ++(7,-7);
\draw[ultra thick,->] (0,13) -- ++(2,-2) -- ++(0,-2) -- ++(2,-2) -- ++(0,-2) -- ++(7,-7);
\end{tikzpicture}
\caption{Inequalities from puzzles.}
\lbl{hex-inequalities}
\end{figure}
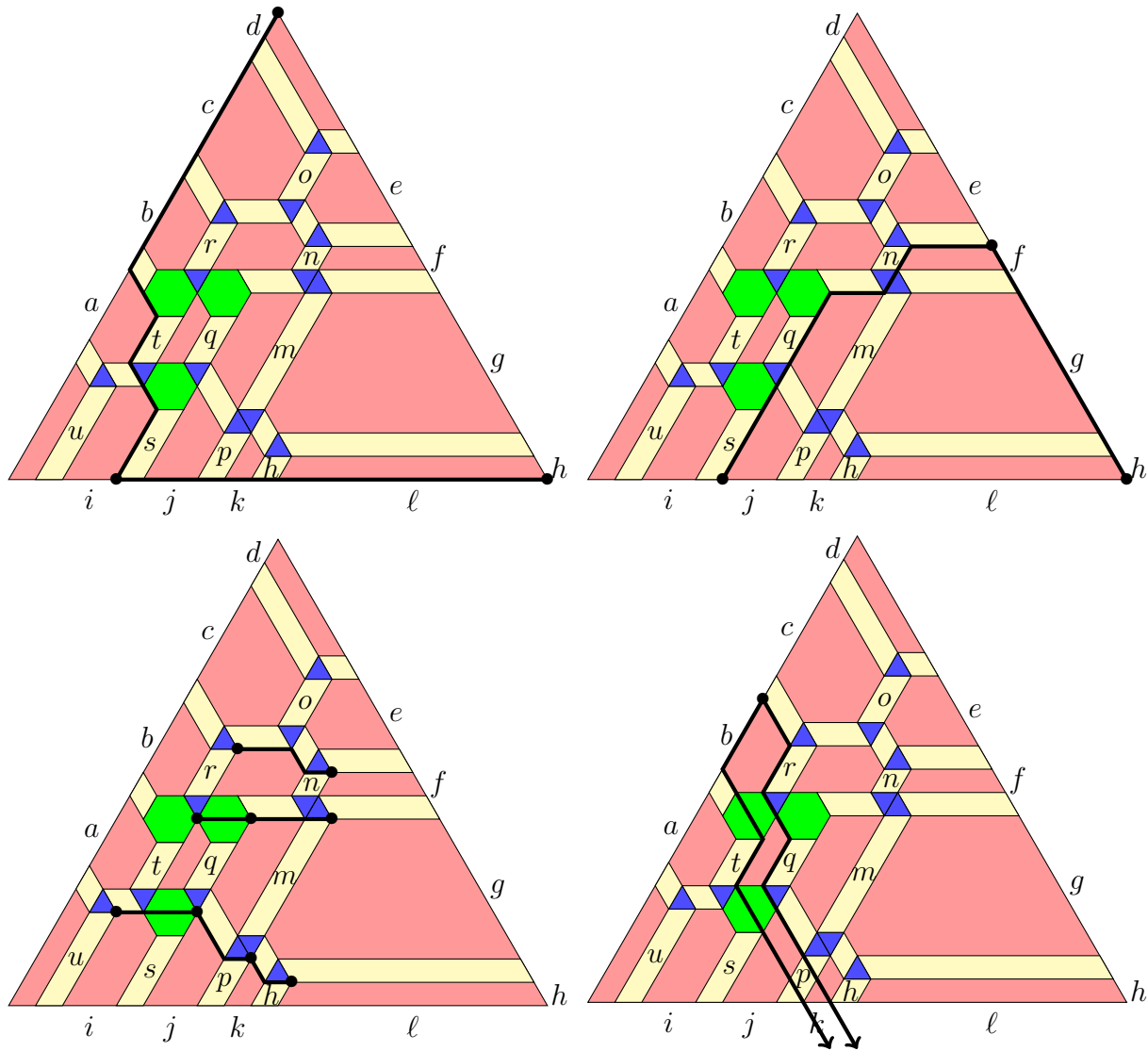

Finally, uncircle the circled numbers in~$T$.
Since circled numbers either share boxes with uncircled numbers
or occur in the inner corners of $\la$,
what we get is a set-valued tableau of shape $\nu/\la^-$,
where $\la^-$ is $\la$ with some inner corners removed.
This concludes one direction of the bijection.

Reversing the bijection is straightforward.
First, we reverse the last step.
Let $T'$ be a set-valued tableau of shape $\nu/\la^-$ and content $\mu$,
where $\la^-$ is $\la$ with some of its inner corners removed.
Circle all the numbers in boxes corresponding to inner corners of $\la$
and all but the smallest number in each of the boxes corresponding to $\nu/\la$.
This tableau with circles is in fact $T$ as described in the middle of the bijection above.
Indeed,
as $T'$ is ballot,
the numbers appearing in Row~$i$ are all at most~$i$.
Also, if we were to get two $\circled i$ in some row,
the right $\circled i$ is sharing its box with a smaller number,
so this row is not weakly increasing from left to right,
a contradiction to the fact that $T'$ is semistandard.

It remains to assemble the puzzle from the tableau~$T$ by reversing the first half of the bijection.
From bottom to top, add in beams of rhombi of the correct height based on the multiplicities of numbers in the tableau,
place an upright $1$-triangle or hexagon on top of each beam depending on the existence of a corresponding circled number,
and join these together using rhombi and upside-down $1$-triangles in the only way possible.
Repeat with the next set of beams and such.
Fill the remaining region with $0$-triangles.
This construction works,
and no tiles need to overlap or extend beyond the boundary,
exactly because the inequalities we derived above
are satisfied if they came from such a tableau.
Checking the details is routine and therefore omitted.

\subsection{Proof of Theorem~\ref{thm:hex-tab-oplus}}

We prove Theorem~\ref{thm:hex-tab-oplus} by establishing a bijection between these tableaux and $\hex$-puzzles.
This bijection is extremely similar to the bijection in the previous proof.
We follow the same outline and use the same running examples.

Given a $\hex$-puzzle, in the boxes of a diagram of shape $\mu\oplus\la$, fill out according to the following schematic plan
\begin{center}
\begin{ytableau}
\none & \none & \none & \none & 1d & 1c & 1b & 1a\\
\none & \none & \none & \none & 2d & 2c & 2b\\
\none & \none & \none & \none & 3d & 3c\\
\none & \none & \none & \none & 4d\\
1u & 2t & 3r & 4o\\
2s & 3q & 4n\\
3p & 4m\\
4h
\end{ytableau}
\end{center}
where, as before,
a number $x$ followed by a letter $y$ means to write $x$ in $y$ adjacent boxes.
If $y'>y$, write an additional $x$ in the \emph{next} box, in its own space.
Circle such a number.
Note that every box has exactly one number, which may or may not be circled.
Call this tableau~$T$.

\begin{example}
Applying the bijection described to the puzzle results in the following tableau.
\begin{center}
\begin{ytableau}
\none & \none & \none & \none & \none & \none & \none & \none & \none & \none & \none & 1 & 1 & 1 & 1 & 1 & 1 & 1 & 1 & 1 & 1 & 1\\
\none & \none & \none & \none & \none & \none & \none & \none & \none & \none & \none & 2 & 2 & 2 & 2 & 2 & 2 & 2 & 2\\
\none & \none & \none & \none & \none & \none & \none & \none & \none & \none & \none & 3 & 3 & 3 & 3 & 3\\
\none & \none & \none & \none & \none & \none & \none & \none & \none & \none & \none & 4\\
1 & 1 & 1 & 1 & 2 & 2 & \circled2 & 3 & 3 & 4 & 4\\
2 & 2 & 2 & \circled2 & 3 & 3 & \circled3 & 4\\
3 & 3 & 4 & 4 & 4 & 4 & 4\\
4\\
\end{ytableau}
\end{center}
\end{example}

We read off \emph{exactly} the same equalities and inequalities from \fig{hex-inequalities}.
However, we interpret them differently.

\textbf{Content.}
The top left picture shows that $\nu_2-\la_2=s+t$,
leading to $\content(T)=\nu$ where $\circled i$ is ignored.

\textbf{Shape.}
The top right picture shows that $\mu_2=s'+q'+n'$,
showing that $\circled i$ shall occupy its own box.

\textbf{Semistandard.}
The lower left picture shows that
$u'+t'\geq s'+q'\geq p'+m'$.
This directly translates to the semistandard condition of~$T$, where $\circled i$ is treated as $i$.

\textbf{Ballot.}
The lower right picture shows the final type of inequalities,
whose interpretation is still slightly more complicated.
Following the notation from the previous proof,
we get $b+t\geq_s r+q$ as one of these inequalities.
Let us see how this kind of inequalities interact with the ballot condition.
Let $w$ be an initial segment of $\row(T)$.
As an example, let us compare the number of $2$s and $3$s.
We may as well extend $w$ with some more $3$s without adding $2$s.
For example, suppose $w$ ends between the $2$s and $3$s of Row~$2$.
There are at least as many (uncircled) $2$s as (uncircled) $3$s in $w$
if and only if $b+t\geq r+q$.
If there is a $\circled2$ between the $2$s and $3$s of Row~$2$,
the incremented erasure of $w$
would have an extra~$3$.
Therefore we must have $b+t\geq_s r+q$.
Other requirements of the ballot condition all amount to inequalities of this type.

This establishes one direction of the bijection.
As before, reversing the bijection and proving correctness is straight-forward, so we omit the details.

\subsection{Bijection between puzzles and tableaux}
Rather than repeat similar proofs over and over,
we present in table form
the inequalities that can be read off from puzzles
and their corresponding interpretations in both skew and $\oplus$ tableaux rules.

\renewcommand{\g}[1]{\geq_{#1}}
\newcommand{\ml}[1]{\begin{tabular}{@{}l@{}}#1\end{tabular}}
\newcommand{\as}[1]{$\circled i\mapsto #1$}

For $\tri$, like for $\hex$, we let $x'=x+1$ if the added tile is above the $x$-beam; otherwise $x'=x$.
For $\K$ and $\hexR$, replace ``above'' in the definition above with ``below.''
Consequently, $u',s',p',h'$ are undefined for $\K$ and $\hexR$.%
\footnote{The $\K$ and $\hexR$ are ``upside down'' in the sense that they replace the upside down $1$-triangle $\star$.
Heuristically, since there are fewer opportunities to use these tiles,
their corresponding set-valued tableaux have no option to fill a larger shape
and circle tableaux have no $\circled i$ in Row~$i$.
The rules in the tables where this manifests itself are marked with ($\star$).}
As before, $x\g zy$ is a shorthand for $x\geq y+z'-z$.

We first re-describe $\hex$ rules in table form to help orient the reader.
\medskip

\begin{tabular}{l|l|l}
\hex & $\nu/\la$ & $\mu\oplus\la$

\\ \hline

{$\begin{aligned}
\nu_1-\la_1 &= u\\
\nu_2-\la_2 &= s+t\\
\nu_3-\la_3 &= p+q+r\\
\nu_4-\la_4 &= h+m+n+o\\
\end{aligned}$}
& \ml{Shape:\\$\circled i$ takes no space \\set-valued}
& \ml{Content:\\ignore $\circled i$}

\\ \hline

{$\begin{aligned}
\mu_1 &= u'+t'+r'+o'\\
\mu_2 &= s'+q'+n'\\
\mu_3 &= p'+m'\\
\mu_4 &= h'\\
\end{aligned}$}
& \ml{Content:\\\as i}
& \ml{Shape:\\$\circled i$ takes a box}

\\ \hline

{$\begin{aligned}
&u'\geq s'\geq p'\geq h'\\
&u'+t'\geq s'+q'\geq p'+m'\\
&u'+t'+r'\geq s'+q'+n'\\
\end{aligned}$}
& \ml{Ballot:\\\as i}
& \ml{Semistandard:\\\as i}

\\ \hline
\begin{tabular}{c}
$a\g ut$ \qquad $b\g tr$ \qquad $c\g ro$ \qquad $d\g o0$\\
$b+t\g sr+q$ \quad $c+r\g qo+n$ \quad $d+o\g n0$\\
$c+r+q\g po+n+m \quad d+o+n\g m0$\\
$d+o+n+m\g h0$\\
\end{tabular}
& \ml{Semistandard:\\$\circled i$ in previous box\\shape becomes $\nu/\la^-$}
& \ml{Ballot:\\ keep only last $\circled i$\\\as{i+1}}
\end{tabular}

\bigskip
The inequalities for $\hexR$ are very similar to those for $\hex$.
The main difference is seen in the last rows of the tables.
Consider the semistandard condition of the skew rule.
While the inequality $a\g ut$ dictates that $\circled1$ in Row~$1$ is to be written in the box \emph{before} the $1$s corresponding to~$u$,
the inequality $a\g tt$ instead dictates that $\circled1$ in Row~$2$ is to be written in the box \emph{after} the $1$s corresponding to~$t$.
Similarly, for the $\oplus$ rule's ballot condition, the erasure is not incremented.
The other difference is marked with ($\star$) due to $\hexR$ being upside down.
\medskip

\begin{tabular}{l|l|l}
\hexR & $\nu/\la$ & $\mu\oplus\la$

\\ \hline

{$\begin{aligned}
\nu_1-\la_1 &= u\\
\nu_2-\la_2 &= s+t\\
\nu_3-\la_3 &= p+q+r\\
\nu_4-\la_4 &= h+m+n+o\\
\end{aligned}$}
& \ml{Shape:\\$\circled i$ takes no space \\set-valued}
& \ml{Content:\\ignore $\circled i$}

\\ \hline

{$\begin{aligned}
\mu_1 &= u+t'+r'+o'\\
\mu_2 &= s+q'+n'\\
\mu_3 &= p+m'\\
\mu_4 &= h\\
\end{aligned}$}
& \ml{Content:\\\as i}
& \ml{Shape:\\$\circled i$ takes a box\\no $\circled i$ in Row $i$ ($\star$)}

\\ \hline

{$\begin{aligned}
&u\geq s\geq p\geq h\\
&u+t'\geq s+q'\geq p+m'\\
&u+t'+r'\geq s+q'+n'\\
\end{aligned}$}
& \ml{Ballot:\\\as i}
& \ml{Semistandard:\\\as i}

\\ \hline
\begin{tabular}{c}
$a\g tt$ \qquad $b\g rr$ \qquad $c\g oo$\\
$b+t\g qr+q$ \quad $c+r\g no+n$\\
$c+r+q\g mo+n+m$\\
\end{tabular}
& \ml{Semistandard:\\$\circled i$ in next box\\stay within shape ($\star$)}
& \ml{Ballot:\\ keep only last $\circled i$\\\as i}
\end{tabular}

As compared to $\hex$, the inequalities for $\K$ look quite different on the surface.
However, it turns out we are essentially swapping the skew and $\oplus$ rules with each other.
Indeed, the only other difference is that $\K$, being upside down, is less frequently usable,
as denoted by ($\star$) in two places.
\medskip

\begin{tabular}{l|l|l}
\K & $\nu/\la$ & $\mu\oplus\la$

\\ \hline

{$\begin{aligned}
\nu_1-\la_1 &= u\\
\nu_2-\la_2 &= s+t'\\
\nu_3-\la_3 &= p+q'+r'\\
\nu_4-\la_4 &= h+m'+n'+o'\\
\end{aligned}$}
& \ml{Shape:\\$\circled i$ takes a box\\no $\circled i$ in Row $i$ ($\star$)}
& \ml{Content:\\\as{i}}

\\ \hline

{$\begin{aligned}
\mu_1 &= u+t+r+o\\
\mu_2 &= s+q+n\\
\mu_3 &= p+m\\
\mu_4 &= h\\
\end{aligned}$}
& \ml{Content:\\ignore $\circled i$}
& \ml{Shape:\\$\circled i$ takes no space \\set-valued}

\\ \hline

{$\begin{aligned}
&u\g ts\g qp\g mh\\
&u+t\g rs+q\g np+m\\
&u+t+r\g os+q+n\\
\end{aligned}$}
& \ml{Ballot:\\ keep only last $\circled i$\\\as{i+1}}
& \ml{Semistandard:\\$\circled i$ in previous box\\stay within shape ($\star$)}

\\ \hline
\begin{tabular}{c}
$a\geq t'$ \qquad $b\geq r'$ \qquad $c\geq o'$\\
 $b+t'\geq r'+q'$ \quad $c+r'\geq o'+n'$\\
$c+r'+q'\geq o'+n'+m'$
\end{tabular}
& \ml{Semistandard:\\\as i}
& \ml{Ballot:\\\as i}
\end{tabular}

\bigskip
The close relation between $\tri$ and $\K$ is similar to that of between $\hex$ and $\hexR$.
Indeed, one difference of $\tri$ compared to $\K$ is that its erasure is not incremented and $\circled i$ goes in the next box,
just like $\hexR$.
On the other hand, the other difference is that $\tri$ does not have ($\star$) restrictions,%
\footnote{So, in the $\oplus$ rule, $\circled i$ can be written in the next box,
even protruding beyond the shape $\mu$.
However, if $\mu_2=\mu_3$, say, the inequalities $s\g pp$ and $s+q\g mp+m$ prohibit $\circled 3$ and $\circled 4$, respectively, from protruding in Row~$3$.
As such, $\mu^+$ is $\mu$ with some \emph{outer corners} added.}
like $\hex$.
The lack of perfect symmetry is somewhat puzzling.
\medskip

\begin{tabular}{l|l|l}
\tri & $\nu/\la$ & $\mu\oplus\la$

\\ \hline

{$\begin{aligned}
\nu_1-\la_1 &= u'\\
\nu_2-\la_2 &= s'+t'\\
\nu_3-\la_3 &= p'+q'+r'\\
\nu_4-\la_4 &= h'+m'+n'+o'\\
\end{aligned}$}
& \ml{Shape:\\\rlap{$\circled i$ takes a box}\hphantom{no $\circled i$ in Row $i$ ($\star$)}}
& \ml{Content:\\\as{i}}

\\ \hline

{$\begin{aligned}
\mu_1 &= u+t+r+o\\
\mu_2 &= s+q+n\\
\mu_3 &= p+m\\
\mu_4 &= h\\
\end{aligned}$}
& \ml{Content:\\ignore $\circled i$}
& \ml{Shape:\\$\circled i$ takes no space \\set-valued}

\\ \hline

{$\begin{aligned}
&u\g ss\g pp\g hh\\
&u+t\g qs+q\g mp+m\\
&u+t+r\g ns+q+n\\
\end{aligned}$}
& \ml{Ballot:\\ keep only last $\circled i$\\\as i}
& \ml{Semistandard:\\$\circled i$ in next box\\shape becomes $\mu^+$}

\\ \hline
\begin{tabular}{c}
$a\geq t'$ \qquad $b\geq r'$ \qquad $c\geq o'$\\
 $b+t'\geq r'+q'$ \quad $c+r'\geq o'+n'$\\
$c+r'+q'\geq o'+n'+m'$
\end{tabular}
& \ml{Semistandard:\\\as i}
& \ml{Ballot:\\\as i}
\end{tabular}


\subsection{Correspondence to coefficients}
In the previous section, we presented in table form the relevant parts of the bijection between
the four puzzle rules given in Section~\ref{ss:K-puzzles} and
the eight tableau rules given in Section~\ref{ss:K-tableaux}.
What remains is to relate these to the coefficients defined in Section~\ref{ss:K-coefficients}.

Buch~\cite{Buch} proved Theorem~\ref{thm:K-tab-oplus}, establishing that the $\K$ rules count~$\c$.
We proved above that the $\hex$ rules count~$\d$.

\subsection{Proof of Theorem~\ref{thm:tri-tab-oplus}}
By definition, we have
\[G_\mu \cdot G_1=\sum_{\mu'}(-1)^{\abs{\mu'}-\abs\mu-1} c_{\mu1}^{\mu'}G_{\mu'}.\]
By Theorem~\ref{thm:K-tab-oplus}, the coefficient $c_{\mu1}^{\mu'}$ is $1$
if $\mu'$ is $\mu$ with a positive number of outer corners added,%
\footnote{Consider the shape $1\oplus\mu$.
The numbers filled in the lower box corresponds to the rows of $\mu'/\mu$.}
and $0$ otherwise.
So
\begin{align*}
G_\la\cdot (G_\mu\cdot G_1)
&=G_\la\sum_{\mu'}(-1)^{\abs{\mu'}-\abs\mu-1}G_{\mu'}\\
&=\sum_{\mu'}(-1)^{\abs{\mu'}-\abs\mu-1} \sum_\nu (-1)^{\abs\nu-\abs\la-\abs{\mu'}}c_{\la\mu'}^\nu G_\nu\\
&=-\sum_{\nu,\mu'}\sign c_{\la\mu'}^\nu G_\nu,
\end{align*}
where $\mu'$ runs over $\mu$ with a positive number of outer corners added.
By definition, we have
\begin{align*}
\sum_\nu\sign\cc G_\nu
&= G_\la\cdot G_\mu\cdot (1-G_1)\\
&=\sum_\nu\sign\c G_\nu + \sum_{\nu,\mu'}\sign c_{\la\mu'}^\nu G_\nu,
\end{align*}
so \[ \cc=\c+\sum_{\mu'}c_{\la\mu'}^\nu. \]
By Theorem~\ref{thm:K-tab-oplus},
$\cc$ is the number of semistandard ballot set-valued tableaux of shape $\mu^+\oplus\la$ and content $\nu$,
where $\mu^+$ is either $\mu$ or $\mu$ with a positive number of outer corners added, as desired.

\subsection{Proof of Theorem~\ref{thm:hexR-tiling}}
By Theorem~\ref{thm:hex-tiling}, it suffices to show a bijection between $\hexR$-puzzles with boundary $\DD$
and $\hex$-puzzles with boundary $\Delta_{\la'\mu'}^{\nu'}$.
The bijection is simple: mirror the puzzle across a vertical line and swap the $0$ and $1$ labels.
This is clearly an involution.
Each of the original puzzle pieces is mapped to a valid puzzle piece.
The $\hex$ and $\hexR$ pieces are mapped to each other.
The boundary is mapped from $\DD$ to $\Delta_{\mu'\la'}^{\nu'}$.%
\footnote{Indeed, recall that the binary string of a partition $\la$ corresponds to the boundary of the diagram of $\la$.
Reversing the string rotates (the boundary of) the diagram by $180^\circ$.
Swapping $0$ and $1$ in the string flips the diagram across the line $x=y$.
Composing these two transformations flips the diagram across the line $x+y=0$.}
Finally, by definition, $\d=d_{\mu\la}^\nu$, so we are done.

\section{Final remarks} \label{sec:4}
\subsection{} \label{ss:protrude}
Consider the example
$\la=(2,1)$, $\mu=(4,2)$, and $\nu=(4,3,1)$.
The skew tableau
\begin{center}
\begin{ytableau}
*(gray) & *(gray) \circled 1& 1 & 1 \\
*(gray) & 2 & 2 \\
1\\
\end{ytableau}
\end{center}
corresponds to the $\hex$-tiling
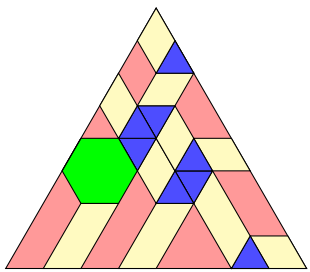
\begin{center}
\begin{tikzpicture}[scale=0.5]
\coords
\draw[fill=\zerocolour] (0,0) -- (8,0) -- (0,8) -- cycle;
\draw[fill=\rhombuscolour] (6,0) rectangle +(1,0);
\draw[fill=\rhombuscolour] (6,0) ++(0,0) ++(0,1) -- ++(1,-1) -- (8,0) -- ++(-1,1) -- cycle;
\draw[fill=\onecolour] (6,0) -- +(1,0) -- +(0,1) -- cycle;
\draw[fill=\rhombuscolour] (6,0) -- ++(0,1) -- ++(-2,2) -- ++(0,-1) -- cycle;
\draw[fill=\onecolour] (4,2) -- ++(0,1) -- ++(-1,0) -- cycle;
\draw[fill=\rhombuscolour] (3,0) rectangle +(1,2);
\draw[fill=\rhombuscolour] (3,2) ++(0,0) ++(0,1) -- ++(1,-1) -- (4,2) -- ++(-1,1) -- cycle;
\draw[fill=\onecolour] (3,2) -- +(1,0) -- +(0,1) -- cycle;
\draw[fill=\rhombuscolour] (3,2) -- ++(0,1) -- ++(-1,1) -- ++(0,-1) -- cycle;
\draw[fill=\onecolour] (2,3) -- ++(0,1) -- ++(-1,0) -- cycle;
\draw[fill=\rhombuscolour] (1,0) rectangle +(1,2);
\draw[fill=\rhombuscolour] (1,3) ++(0,0) ++(0,1) -- ++(1,-1) -- (2,3) -- ++(-1,1) -- cycle;
\draw[fill=green] (1,2) -- +(1,0) -- +(1,1) -- +(0,2) -- +(-1,2) -- +(-1,1) -- cycle;
\draw[fill=\rhombuscolour] (0,3) -- ++(0,1) -- ++(-0,0) -- ++(0,-1) -- cycle;
\draw[fill=\rhombuscolour] (3,3) rectangle +(1,0);
\draw[fill=\rhombuscolour] (3,3) ++(0,0) ++(0,1) -- ++(1,-1) -- (5,3) -- ++(-1,1) -- cycle;
\draw[fill=\onecolour] (3,3) -- +(1,0) -- +(0,1) -- cycle;
\draw[fill=\rhombuscolour] (3,3) -- ++(0,1) -- ++(-1,1) -- ++(0,-1) -- cycle;
\draw[fill=\onecolour] (2,4) -- ++(0,1) -- ++(-1,0) -- cycle;
\draw[fill=\rhombuscolour] (1,4) rectangle +(1,0);
\draw[fill=\rhombuscolour] (1,4) ++(0,0) ++(0,1) -- ++(1,-1) -- (2,4) -- ++(-1,1) -- cycle;
\draw[fill=\onecolour] (1,4) -- +(1,0) -- +(0,1) -- cycle;
\draw[fill=\rhombuscolour] (1,4) -- ++(0,1) -- ++(-1,1) -- ++(0,-1) -- cycle;
\draw[fill=\rhombuscolour] (1,5) rectangle +(1,1);
\draw[fill=\rhombuscolour] (1,6) ++(0,0) ++(0,1) -- ++(1,-1) -- (2,6) -- ++(-1,1) -- cycle;
\draw[fill=\onecolour] (1,6) -- +(1,0) -- +(0,1) -- cycle;
\draw[fill=\rhombuscolour] (1,6) -- ++(0,1) -- ++(-1,1) -- ++(0,-1) -- cycle;
\end{tikzpicture}
\end{center}
Since the shapes all fit in a $4\times3$ box, one might think $n=7$ is sufficient side length for a puzzle.
However, the $\hex$ piece will protrude to the left of the puzzle with side length~$7$.
In Theorem~\ref{thm:hex-tiling} (and analogously in Theorem~\ref{thm:hexR-tiling}),
we dealt with this issue by increasing the puzzle size length by one.

Another way to solve this issue is to add an additional trapezoid piece
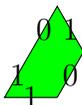
\begin{center}
\begin{tikzpicture}[scale=0.7]
\begin{scope}[shift={(7,0)}]
   \draw[fill=green] (0,0) \r1\b2\r4\b4\b0;
\end{scope}
\end{tikzpicture}
\end{center}
as if to allow the hexagonal tile to protrude to the left.
(By the way things are set up, the hexagon never needs to protrude to the right or below.)
However, we do not want this piece used elsewhere.
So we must make some more modifications.
Here is the complete set of tiles.

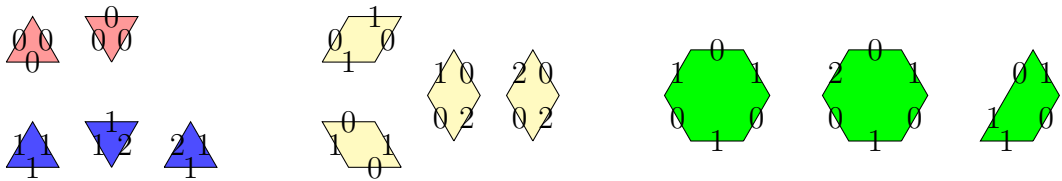
\begin{figure}[hbtp]
\begin{tikzpicture}[scale=0.7]
\draw[fill=\zerocolour] (0,2)
   \r0\r2\r4;
\draw[fill=\zerocolour] (2,2)
   \r1\r3\r5;

\draw[fill=\onecolour] (0,0)
   \b0\b2\b4;
\draw[fill=\onecolour] (3,0)
   \b0\b2\bb4;
\draw[fill=\onecolour] (2,0)
   \bb1\b3\b5;

\begin{scope}[shift={(6,0)}]
\draw[fill=\rhombuscolour] (0,2)
   \b0\r1\b3\r4;
\draw[fill=\rhombuscolour] (1.5,0)
   \b2\r3\b5\r0;
\draw[fill=\rhombuscolour] (2.5,0.5)
   \bb1\r2\b4\r5;
\draw[fill=\rhombuscolour] (4,0.5)
   \bb1\r2\bb4\r5;
\end{scope}

\begin{scope}[shift={(14,0.5)}]
   \draw[fill=green] (0,0) \r1\b2\r3\b4\r5\b0;
   \draw[fill=green] (3,0) \r1\b2\r3\bb4\r5\b0;
   \draw[fill=green] (5.5,0) \r1\b2\r4\b4\b0;
\end{scope}
\end{tikzpicture}
\end{figure}

Consider the northeast--southwest slanting $1$ edges.
A $1$-edge on the bottom-right side of pieces are now labelled with $2$,
so the new trapezoid piece cannot be used except at the left boundary.
An old piece with a $1$-edge on its top-left side
must be duplicated,
with a version for use at the left boundary and another for use in the interior.

\subsection{} \label{ss:genomic}
Theorem~\ref{thm:K-tab-skew} gives a skew tableau rule for calculating the $K$-theoretic Littlewood--Richardson coefficients $\c$
using right circle tableaux.
Pechenik and Yong give the same rule using \emph{genomic} tableaux (see \cite{genomic} for definitions).

\begin{thm}[$\K$, Pechenik--Yong~\cite{genomic}, $K$-theory, skew version] \label{thm:K-genomic} 
The coefficient $\c$ is the number of semistandard ballot \uline{genomic} tableaux of shape $\nu/\la$ and content $\mu$.
\end{thm}

These two rules are virtually identical, as there is a simple bijection between
right circle tableaux and genomic tableaux.
Indeed,
let a semistandard ballot right circle tableau of shape $\nu/\la$ and content $\mu$ be given.
By semistandardness, the boxes filled with $i$ and $\circled i$ form a horizontal strip.
From left to right, rewrite these as $i_1$, $i_2$, $i_3$, and so on.
Whenever $\circled i$ is encountered, the \emph{next} subscript used is the same as this subscript.
By ballotness, the rightmost $i$ in the tableau is not circled, so this rule is well-formed.
It is easy to see that this yields a semistandard genomic tableau of the same shape and content.
One can also check that the tableau is ballot.

Conversely,
given a semistandard ballot genomic tableau,
the boxes filled with $i_j$ for a fixed $i$ form a horizontal strip.
From left to right, circle an entry if its subscript is the same as the next one.
Erase all subscripts.
The correctness of this bijection is straightforward and left as exercise to the reader.

\begin{example}
The structure constant $c_{(2,1),(2,1)}^{(4,2,1)}$ is computed by the circle tableaux
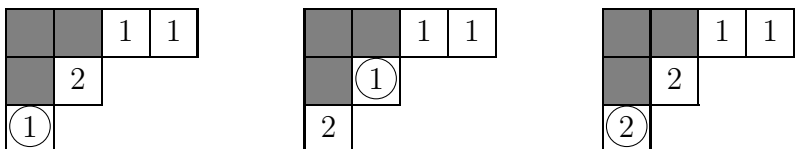
\begin{center}
\begin{ytableau}
*(gray) & *(gray) & 1 & 1 \\
*(gray) & 2 \\
\circled1\\
\end{ytableau}\hspace{.5in}
\begin{ytableau}
*(gray) & *(gray) & 1 & 1 \\
*(gray) & \circled1 \\
2\\
\end{ytableau}\hspace{.5in}
\begin{ytableau}
*(gray) & *(gray) & 1 & 1 \\
*(gray) & 2 \\
\circled2\\
\end{ytableau}
\end{center}
and by the corresponding genomic tableaux
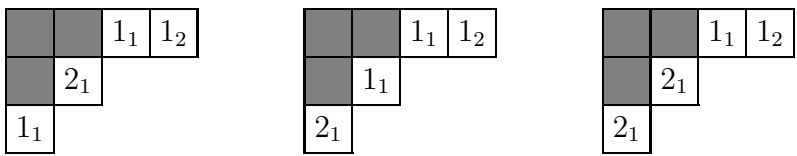
\begin{center}
\begin{ytableau}
*(gray) & *(gray) & 1_1 & 1_2 \\
*(gray) & 2_1 \\
1_1\\
\end{ytableau}\hspace{.5in}
\begin{ytableau}
*(gray) & *(gray) & 1_1 & 1_2 \\
*(gray) & 1_1 \\
2_1\\
\end{ytableau}\hspace{.5in}
\begin{ytableau}
*(gray) & *(gray) & 1_1 & 1_2 \\
*(gray) & 2_1 \\
2_1\\
\end{ytableau}
\end{center}
\end{example}

\newpage
\subsection*{Acknowledgements}
We are grateful to
Allen Knutson and
Joel Lewis
for helpful conversations.
The first-named author is partially supported by NSF grant DMS-1351590 and Sloan Fellowship.
The second-named author is partially supported by NSF RTG grant NSF/DMS-1148634.

\bibliography{ktiles}
\bibliographystyle{htam} 

\end{document}